\documentclass[leqno,12pt,draft]{article}
\usepackage{amssymb}
\begin{document}

{\sc Krzysztof Zajkowski }
\begin{center}
{\Large Covariance matrices of self-affine measures}
\end{center}
\begin{abstract}
In this paper we derive a formula for a covariance matrix
of any self-affine measure, i.e. a probability measure
$\mu$ satisfying
$$ \mu = \sum_{k=1}^l p_k\mu \circ S_k^{-1}\;, $$
where $\{S_k({\bf x})=A_k {\bf x} + {\bf b}_k\}_{1\leq k \leq l}$ 
is a family
of affine contractive maps and $\{p_k\}_{1\leq k \leq l}$
is a set of probability weights. In particular if for every $k$,
$A_k=A$ then the formula will be have the following form
$$ D^2X=[I\otimes I-A\otimes A]^{-1}\mathcal{D}^2\mathcal{B},$$
where $D^2X$ denote the covariance matrix of the measure $\mu$
and $\mathcal{D}^2\mathcal{B}$ denote a covariance matrix of 
a discret random variable $\mathcal{B}$
with values ${\bf b}_k$, and corresponding probabilites $p_k$. 
\end{abstract}


\setcounter{section}{0}
\section{Introduction}

We will consider a invariant probability measure on ${\Bbb R}^d$
\begin{equation}
\label{dfm}
\mu = \sum_{k=1}^l p_k\mu \circ S_k^{-1}\;,
\end{equation}
where $\{S_k\}_{1\leq k \leq l}$ is a family of contractive maps
and $\{p_k\}_{1\leq k \leq l}$ is a set of probability weights. 
It is often assumed
that the maps are similitudes. We will make a more general
assumption that the maps are affine contractions on ${\Bbb R}^d$;
i.e. $S_k({\bf x})=A_k{\bf x}+{\bf b}_k$ 
and the operator norm
$\Vert A_k \Vert< 1$ for all $k$.

The definition (\ref{dfm}) was introduced by 
Hutchinson [H]. But an example of such measures has been studied 
for a long time in the context of Bernoulli convolution,
i.e. the example of an invariant measure on the real line
\begin{equation}
\mu=\frac{1}{2}(\mu\circ S^{-1}_1+\mu\circ S^{-1}_2),
\end{equation}
where $S_1(x)=\beta(x+1)$, $S_2(x)=\beta(x-1)$ for 
$\beta\in(0,1)$.
It remains difficult open problem to characterize the set
of $\beta$ for which $\mu$ is absolutely continuous.
An another example has been studied in great detail
in wavelet theory in connection with the dilatation
equation
\begin{equation}
f(x)=\sum_{k=1}^lc_k f(2x-(k-1)).
\end{equation}
The function $f$ can be considered as the density function 
of the corresponding absolutely continuous self-affine
measure $\mu$ for $S_k(x)=\frac{1}{2}(x+(k-1))$ and
$p_k=\frac{1}{2}c_k$. In wavelet theory the $c_k$ may be
negative but $\sum c_k$ must be $2$.
The invariant measures have a natural connection with
fractal geometry [F],
because their supports are compact invariant sets
\begin{equation}
K=\bigcup_{k=1}^lS_k(K).
\end{equation}
These measures arise also in another areas of mathematics.
 
A fundamental method for studing these measures is
the Fourier transform [S]. 
Our goal is to show that we can use another 
probabilistic tools, not only characteristic function,
to investigate the self-affine measures. We derive a formula
for  covariance matrices and give an example of an investigation
of measures on Sierpinski triangle.

\section{Covariance matrices of self-affine measures}
The invariant measure $\mu$ satisfies the following identity
\begin{equation}\label{int}
\int_{{\Bbb R}^d}f d\mu=\sum_{k=1}^l p_k\int_{{\Bbb R}^d}
f\circ S_k d\mu,
\end{equation}
where $f$ is any continuous function on ${\Bbb R}^d$ [B]. 
We will denote by $X$ some $d$-dimensional random variable 
with respect to
the probability distribution $\mu$. 
We apply the identity (\ref{int}) to the coordinate
functions $e_i^\star({\bf x})=x_i$ of the point
${\bf x}=(x_i)_{1\leq i \leq d}\in {\Bbb R}^d$.
By the above a vector of expected values $EX$ will be equal

\begin{eqnarray}
EX & = & (\int_{{\Bbb R}^d}e_i^\star({\bf x}) d\mu)_{1\leq i \leq d}
=(\sum_{k=1}^l p_k \int_{{\Bbb R}^d}
e_i^\star(A_k{\bf x}+{\bf b}_k)d\mu)
_{1\leq i \leq d}\nonumber\\
\; & = & \sum_{k=1}^l p_k E(A_kX+{\bf b}_k) 
= \sum_{k=1}^lp_kA_kEX + \sum_{k=1}^l p_k {\bf b}_k
\end{eqnarray}
The above relation gives
\begin{equation}
\label{ex}
[I-\sum_{k=1}^l p_k A_k]EX = \sum_{k=1}^l p_k {\bf b}_k,
\end{equation}
where $I$ is identity matrix on ${\Bbb R}^d$.
The sum $\sum_{k=1}^l p_k {\bf b}_k$ is a vector of 
expected values of
a $d$-dimensional random variable $\mathcal{B}$ with values ${\bf b}_k$,
and corresponding probabilities $p_k$. 
Since  $\Vert A_k \Vert < 1$ for all $k$ then $\Vert 
\sum_{k=1}^l p_k A_k \Vert <1$. It follows
that 1 is not an eigenvalue
of the operator $\sum_{k=1}^l p_k A_k$. For this reason the
operator $I-\sum_{k=1}^l p_k A_k$ will be invertible.   
We can rewrite (\ref{ex}) as 
\begin{equation}
\label{exval}
EX=[I-\sum_{k=1}^l p_k A_k]^{-1}\cal{EB}.
\end{equation}
This means that the expected value of $X$ linearly depend on
the expected value of $\mathcal{B}$.

Let $X\otimes X=[x_ix_j]_{1\leq i,j \leq d}$ denote 
the second order tensor build from the coordinates.
Using (\ref{int}) a matrix of second order moments
$E(X\otimes X)$ will be equal
\begin{eqnarray}
E(X\otimes X)& = & [\int_{{\Bbb R}^d}
e_i^\star({\bf x})e_j^\star({\bf x}) d\mu]
_{1\leq i,j\leq d}\nonumber\\
\; & = & \sum_{k=1}^l p_k [\int_{{\Bbb R}^d}
e_i^\star(A_k{\bf x}+{\bf b}_k)
e_j^\star(A_k{\bf x}+{\bf b}_k)d\mu]_{1\leq i,j\leq d}\nonumber\\
\; & = & \sum_{k=1}^l p_k E((A_k X + {\bf b}_k)
\otimes(A_kX+{\bf b}_k))\nonumber\\
\; & = & \sum_{k=1}^lp_k[(A_k\otimes A_k)E(X\otimes X)+
{\bf b}_k\otimes A_k EX + A_k EX
\otimes {\bf b}_k+ {\bf b}_k\otimes {\bf b}_k]
\end{eqnarray}
Therefore
\begin{equation}
[I\otimes I-\sum_{k=1}^lp_k(A_k\otimes A_k)]E(X\otimes X)
= \sum_{k=1}^l p_k({\bf b}_k\otimes A_kEX+A_kEX\otimes {\bf b}_k
+{\bf b}_k\otimes {\bf b}_k)
\end{equation}
The operator norm of $A_k\otimes A_k$ is less than $1$ on
${\Bbb R}^d\otimes{\Bbb R}^d$,
so, by the same argument as early, we get that the operator
$I\otimes I-\sum_{k=1}^lp_k(A_k\otimes A_k)$ is invertible
and we obtain
\begin{eqnarray}
E(X\otimes X)& =& [I\otimes I-\sum_{k=1}^lp_k(A_k\otimes A_k)]^{-1}
\times \nonumber\\
\; & \; & {\sum_{k=1}^l p_k({\bf b}_k\otimes A_kEX+A_kEX\otimes 
{\bf b}_k + {\bf b}_k\otimes {\bf b}_k)}.
\end{eqnarray}
Substituting (\ref{exval}) into $D^2X=E(X\otimes X)-EX\otimes EX$
we can obtain a formula for the covariance matrix of $X$.
But in the general case it will be complicated.
This formula takes
a surprising simple form when we assumed that all affine
maps have the same linear part, i.e. all $A_k=A$. Under this 
assumption, using the standard tensor calculus we get
\begin{eqnarray}
D^2X & = & [I\otimes I-A\otimes A]^{-1} 
\lbrace [I\otimes A(I-A)^{-1}](\cal{EB}\otimes\cal{EB})\nonumber\\
\; & \; & +[A(I-A)^{-1}\otimes I](\cal{EB}\otimes\cal{EB})
+\cal{E}(\cal{B}\otimes\cal{B})\} \nonumber\\
\; & \; &-[(I-A)^{-1}\otimes (I-A)^{-1}](\cal{EB}\otimes\cal{EB})
\nonumber\\
\; & = & [I\otimes I - A\otimes A ]^{-1}(\cal{E}
(\cal{B}\otimes\cal{B})-\cal{EB}\otimes\cal{EB})
\end{eqnarray}
Notice now that the expression 
$\cal{E}(\cal{B}\otimes\cal{B})-\cal{EB}\otimes\cal{EB}$
is the covariance matrix of $\cal{B}$. Thus we obtain
\begin{equation}
D^2X=[I\otimes I - A\otimes A]^{-1}\mathcal{D}^2\mathcal{B}.
\end{equation} 
In other words we obtained the following proposition.

PROPOSITION

Assume that $\mu$ is a self-affine measure on ${\Bbb R}^d$
for a family of linear contractions $S_k({\bf x })=A{\bf x}
+{\bf b}_k$, $1\leq k \leq l$.
Let $X$ be some random variable with respect
to the probability ditribution $\mu$.
Then the covariance matrix of the random variable $X$
\begin{equation}
\label{cov}
D^2X=[I\otimes I - A\otimes A]^{-1}\mathcal{D}^2\mathcal{B},
\end{equation}
where $\mathcal{D}^2\mathcal{B}$ denote the covariance
matrix of the random variable $\mathcal{B}$.
 $\square$

{\bf Remark.}
When all affine maps $S_k$ have the same linear part then
not only the expected value of $X$ linearly depend on expected
value of $\cal{B}$ but also covariance matrix of $X$ linearly
depend on the covariance matrix of $\cal{B}$.

If the matrix $A$ is diagonal then diagonal is the matrix
$[I\otimes I - A\otimes A]^{-1}$. Therefore we get the
simple corollary.

COROLLARY

If under the assumptions of Proposition
we assume additionaly that the matrix $A$ is diagonal then 
$X_i=e_i^\star(X)$ and $X_j=e_j^\star(X)$ 
are uncorrelated if and only if uncorrelated 
are $(e_i^\star({\bf b}_k))_{1\leq k\leq l}$ and 
$(e_j^\star({\bf b}_k))_{1\leq k\leq l}$.
$\square$

In other words we have obtained the following law
\begin{eqnarray}
\int_{{\Bbb R}^d}x_ix_jd\mu=\int_{{\Bbb R}^d}x_id\mu
\int_{{\Bbb R}^d}x_jd\mu & iff &
\sum_{k=1}^lp_k e_i^\star({\bf b}_k) e_j^\star ({\bf b}_k)
=\sum_{k=1}^lp_ke_i^\star({\bf b}_k)
\sum_{k=1}^lp_ke_j^\star({\bf b}_k).\nonumber
\end{eqnarray}

\section{Example}

Consider Sierpinski triangle with vertices at the points
$(0,0)$, $(1,0)$ and $(\frac{1}{2},\frac{\sqrt{3}}{2})$.
The Sierpinski triangle is an invariant compact set of
three cotractions on ${\Bbb R}^2$: $S_1({\bf x})=\frac{1}{2}{\bf x}$,
$S_2({\bf x})=\frac{1}{2}{\bf x}+(\frac{1}{2},0)$ and
$S_3({\bf x})=\frac{1}{2}{\bf x}+(\frac{1}{4},\frac{\sqrt{3}}{4})$.
In this case the matrix $A=\frac{1}{2}I$. Let ${\mu}$ denote
an invariant measure for weigths $p_1$, $p_2$ and $p_3$.
The expected value $\mathcal{EB}=(\frac{1}{2}p_2+\frac{1}{4}p_3,
\frac{\sqrt{3}}{4}p_3)$. The matrix $[I-A]^{-1}=2I$. By
the (\ref{exval}) we get $\int_{{\Bbb R}^2}x_1d\mu=p_2+\frac{1}{2}p_3$
and $\int_{{\Bbb R}^2}x_2d\mu=\frac{\sqrt{3}}{2}p_3$.
By the (\ref{cov}) we can obtain terms of the matrix
$\mathcal{D}^2\mathcal{B}$ and $D^2X$. In particular
\begin{equation}
\sum_{k=1}^3p_ke_1^\star({\bf b}_k)e_2^\star({\bf b}_k)
-\sum_{k=1}^3p_k e_1^\star({\bf b}_k)
\sum_{k=1}^3p_k e_2^\star({\bf b}_k) 
=\frac{\sqrt{3}}{16}p_3(p_1-p_2).
\end{equation}
By corollary, if $p_1=p_2$ then
the random variables $X_1$, $X_2$ are uncorrelated and 
\begin{eqnarray}
\int_{{\Bbb R}^2}x_1x_2d\mu & = & \int_{{\Bbb R}^2}x_1d\mu
\int_{{\Bbb R}^2}x_2d\mu\nonumber\\
\; & = & (p_2+\frac{1}{2}p_3)\frac{\sqrt{3}}{2}p_3
=\frac{\sqrt{3}}{4}p_3.
\end{eqnarray}

{\bf References}

[B]  M.F Barnsley, {\it Fractals everywhere}, 2nd edition,
Academic Press, Orlando, FL.1993.

[F]  K.J Falconer, {\it Probabilistic mathods in fractal geometry},
Progress in Probability, Vol 37(1995)3-13.

[H]  J.E Hutchinson, {\it Fractals and self similarity},
Indiana Univ. Math. J. 30(1981)713-747.

[S]  R.S. Strichartz, {\it Self-similar measures and their
Fourier transform III}, Indiana Univ. Math. J. 42(1993)
367-411.

\begin{flushleft}
{\small Krzysztof Zajkowski\\
Institute of Mathematics, University of Bialystok\\ 
Akademicka 2, 15-267 Bialystok, Poland\\
E-mail:{\tt kryza@math.uwb.edu.pl}}
\end{flushleft}

\end{document}